
\magnification=\magstep1
\baselineskip=12pt
\overfullrule=0pt
\input amssym.tex

\font\bigbf = cmbx10 scaled \magstep2
\def\noind{\noindent}

\def\cF{{\cal F}}
\def\cD{{\cal D}}

\def\cO{{\cal O}}
\def\cF{{\cal F}}
\def\cL{{\cal L}}

\def\prf{{\noindent \bf Proof.\quad}}

\def\tworight{{\twoheadrightarrow}}
\def\lright{{\longrightarrow}}

\def\Rarrow{{\Rightarrow}}
\def\hookr{{\hookrightarrow}}

\def\IA{{\Bbb A}}
\def\IC{{\Bbb C}}
\def\IO{{\Bbb O}}
\def\IQ{{\Bbb Q}}
\def\IZ{{\Bbb Z}}

\def\Ker{{\rm Ker}}
\def\Spec{{\rm Spec}}
\def\Im{{\rm Im}}
\def\Tor{{\rm Tor}}

\def\vir{{\rm vir}}
\def\b{{\bullet}}
\def\gr{{\rm gr}}

\font\bigbf = cmbx10 scaled \magstep2

\def\oOmega{{\overline \Omega}}

\def\oE{{\overline E}}

\def\oT{{\overline T}}

\def\vep{{\varepsilon}}

\def\tY{{\tilde Y}}

\def\uH{{{\underline H}}}
\def\uK{{{\underline K}}}

\centerline {\bigbf Virtual fundamental classes via dg-manifolds}

\vskip 1cm

\centerline {\bf Ionu\c t Ciocan-Fontanine and Mikhail Kapranov}

\vskip 2cm

\centerline {\bf Introduction}

\vskip 1cm

\noindent {\bf (0.1)} In many moduli problems in algebraic geometry there is a difference
between the actual dimension of the moduli space and the expected, or virtual
dimension. When this happens, the moduli problem is said to be obstructed. The
actual dimension, at the level of tangent spaces,
is typically the dimension of $H^0$ or $H^1$ of some
coherent sheaf $\cal F$, while the virtual dimension is the  Euler characteristic of $\cal F$.
Over $\IC$, one can often represent the moduli space as a possibly non-transversal
intersection inside an infinite-dimensional ambient space, and by analogy
with the finite-dimension intersection theory [F] one expects a ``virtual fundamental
class'' of the expected dimension, associated to the moduli space. Such classes
were constructed by Behrend and Fantechi [BF], and by Li and Tian [LT],
for the case when the obstruction is simple, or
"perfect" (typically, $\cal F$ has one more cohomology group). In this case
the expected dimension is less or equal than the actual one, and the class
lies in the Chow group of the moduli space. 

\vskip .3cm

\noindent {\bf (0.2)} M. Kontsevich suggested in [K]  that all  such problems can be handled
by working with appropriate derived versions of moduli spaces. Following this suggestion,
the authors developed in [CK1-2] the basic theory of such derived objects, called
dg-manifolds, and constructed the derived versions of Grothendieck's $Quot$ and Hilbert
schemes as well as of Kontsevich's moduli spaces of stable maps.

The goal of the present paper is to define virtual classes in the context of
``simply obstructed'' dg-manifolds. By simply obstructed we mean that the tangent
dg-spaces have cohomology only in degrees 0 and 1. Some of the features
of our approach are similar to those of [BF]. In particular, it is clear
that whenever both approaches are applicable, they give the same result.
On the other hand, the language of dg-manifolds exhibits all the necessary
constructions as analogs of the most standard procedures of  usual algebraic
geometry. In particular, the structure sheaf of a dg-manifold gives rise to the
K-theoretic virtual class, and we prove (Theorem 3.3) that it lies in the right level
of the dimension filtration and gives the homological class after passing to
the quotient. Further, we prove a Riemann-Roch-type result for dg-manifolds
(Theorem 4.5.1) which involves integration over the virtual class. 
In a similar way, applying the Bott-Thomason localization theorem to the structure sheaf
of a dg-manifold with a torus action gives at once the localization theorem
for virtual classes proved by Graber and Pandharipande [GP]. 

The intuitive point of view behind  the
language of  dg-manifolds is that they
provide an algebro-geometric analog of ``deformation to transversal intersection''
which often cannot be achieved within pure algebraic geometry. We prove a result
confirming this intuition in a new way.
Namely, we associate, in Theorem 4.6.4, to each dg-manifold $X$ of our type
a cobordism class of almost complex (smooth)  manifolds.

\vskip .3cm

\noindent {\bf (0.3)} One of the most attractive features of our approach is that
it suggests a definition of the virtual class also in the case when the obstruction
is no longer simple. In this case it is not even clear a priori where the virtual class
should lie, as the expected dimension can well be greater that the actual one
(due to many alternating summands in the Euler characteristic). The language of
dg-manifolds suggests that it should lie in the Chow group (of the expected dimension)
of a certain natural fiber bundle $\Pi$ over the moduli space. To be precise
(see (1.1) below), a dg-manifold $X$ consists of a smooth algebraic variety $X^0$
and a sheaf $\cO^\bullet_X$ of dg-algebras on $X^0$. The role of the moduli space
is played by the subscheme $\pi_0(X)\subset X^0$ which is the spectrum of
$\uH^0(\cO_X^\bullet)$. The fiber bundle $\Pi$ is the spectrum of $\uH^{{\rm even}}(\cO_X^\bullet)$,
the ring of even cohomology, and the odd cohomology gives a coherent sheaf $\cal H$ on it.
The virtual class should lie in the Chow group of $\Pi$ and come from the class
$1-[{\cal H}]$ in its K-theory. Further, since  $\uH^{{\rm even}}(\cO_X^\bullet)$
is graded, $\Pi$ is a cone with apex $\pi_0(X)$, so it has an action of ${\Bbb G}_m$
with fixed locus $\pi_0(X)$. This allows one to localize all the data back to $\pi_0(X)$. 
 This program will be developed
in a future paper. 

\vskip .3cm

\noindent {\bf (0.4)} Here is the outline of the paper. In \S 1 we develop the
formalism of deformation to the normal cone in the context of dg-manifolds. It
allows us to replace, in enumerative arguments, the underlying variety $X^0$ of a dg-manifold
$X$ by the normal cone to $\pi_0(X)$ in $X^0$. In section 2 we introduce the class
of $[0,1]$-manifolds which formalize the concept of a simply obstructed moduli space. 
An important property of such manifolds is that the cohomology $\uH^\bullet(\cO_X^\bullet)$
is bounded, so one can speak about its class in the Grothendieck group of
$\pi_0(X)$. This is exactly the K-theoretic virtual class as defined in \S 3. 
We also introduce, in \S 3, the homological class and compare it to the K-theoretic one.  
In section 4 we give a different definition of the homological virtual class
in terms of the Chern character. This was the initial proposal of Kontsevich [K].
Therefore our paper connects the approaches of [K] and [BF]. 
This equivalence of the two definitions
can be seen as a particular case of a Riemann-Roch theorem for dg-manifolds
which we also prove in \S 4. Finally, \S 5 is devoted to the Bott localization
for dg-manifolds.

\vskip .3cm

\noindent {\bf (0.5)} A large part of this paper was written when the second author
was visiting the University of Minnesota. He would like to thank the University
for hospitality and financial support. In addition, the research of both authors
was partially supported by the NSF. 

\vfill\eject

\baselineskip =15pt

\centerline{\bf 1.\quad Deformation to the normal cone for
dg-manifolds} 

\vskip 1cm

\noind {\bf (1.1)\quad Notation.} We fix a base field $k$ of characteristic 0. 
Recall (see [CK1] for
more background) that a dg-scheme is a dg-ringed space
$X=(X^0, \cO^\bullet_X)$, where $X^0$ is a $k$-scheme
and $\cO^\bullet_X$ is a sheaf of dg-algebras on $X^0$,
situated in degrees $\leq 0$, such that $\cO^0_X = \cO_{X^0}$
and quasicoherent as a module over $\cO_{X^0}$. 
We denote the differential in $\cO^\bullet_X$ by $d$.
Because of the grading condition, $d$ is linear over
$\cO_{X^0}$. Further, $\uH^0(\cO_X^\bullet) = 
\cO_{X^0}/d(\cO^{-1}_X)$ is a quotient of $\cO_{X^0}$,
so $\pi_0(X):=  {\rm Spec} \uH^0(\cO_X^\bullet)$
is a closed subscheme of $X^0$.
A dg-sheaf on a dg-scheme $X$ is a sheaf $\cF^\bullet$ of dg-modules
over $\cO_X^\bullet$ which is quasicoherent over $\cO_{X^0}$. 
A dg-sheaf is called a dg-vector bundle, if it is bounded from above,
and, considered as a sheaf of graded modules over $\cO_X^\bullet$,
 is locally free with finitely many generators in each degree. 

\vskip .2cm

By a dg-manifold we mean a dg-scheme $X$ such that $X^0$ is
a smooth algebraic variety over $k$, and $\cO_X^\bullet$, considered
as a sheaf of graded $\cO_{X^0}$-algebras, is locally free with finitely many
generators in each degree.

\vskip .2cm

 Let $X$ be a dg-manifold. For any dg-vector
bundle $E^\bullet$ on $X$ we denote by $E^\bullet\big |_{{\pi_0}(X)}$ the
restriction of $E^\bullet$ (as a complex of vector bundles
on $X^0$) to $\pi_0(X)\subset X^0$. The restriction 
${\cal O}^\b_X\big |_{{\pi_0}(X)}$
will be denoted by $\IO^\b_X$ or simply $\IO^\b$. This 
is a sheaf of dg-algebras on $\pi_0(X)$ with
$d:\IO^{-1}\to\IO^0$ vanishing. Thus $\IO^{\leq -1}$ is a dg-ideal in $\IO^\b$.
 For any $E^\b$ as above the restriction 
$E^\b\big |_{{\pi_0(X)}}$ is a dg-module over $\IO^\b$.

We denote $\oE ^\b=E^\b\otimes_{{\cal O}^\b_X}
{\cal O}_{{\pi_0}(X)}$ the
restriction of $E^\b$ to ${\pi_0}(X)$ in the sense of
dg-manifolds. This is a complex of vector bundles on
$\pi_0(X)$. It is clear that
$$\oE^\b=E^\b\big|_{{\pi_0}(X)} \biggl/ \IO^{\leq -1} 
\cdot E^\b\big|_{{\pi_0}(X)}.$$

We also denote
$$\omega^\b=\omega^\b_X=\oOmega\,^{1, \b}_X,\ {\bf t}^\b = {\bf t}^\b_X=
\oT\,^\b_X.$$
These are complexes of vector bundles on $\pi_0(X)$ 
situated in degrees $\le 0,\ \ge 0$ respectively, and dual
to each other. Note that ${\bf t}^0=T_{X^0}\big|_{{\pi_0}(X)}$,
while for $n>0$
$${\bf t}^n=\Ker\{\IO^{-n *}\to \bigoplus\limits_{i+j=n\atop i,j>0}
\IO^{-i*}\otimes \IO^{-j*}\}$$
is the space of primitive elements in
$\IO^{-n*}$. In particular, ${\bf t}^1=\IO^{-1*}$.

For a dg-bundle $E^\b$ on $X$ we have the
decomposability filtration $\cal D$ in $E^\b\big|_{{\pi_0}(X)}$
$${\cal D}^n E^\b\big|_{{\pi_0}(X)}=(\IO^{\le-1})^n\cdot 
E^\b\big|_{{\pi_0}(X)}$$

\proclaim (1.1.1) Proposition.  We have
$${\rm gr}^n_\cD \big(E^\b\big|_{{\pi_0}(X)}\big)=
\oE\,^\b\otimes S^n(\omega^{\le-1}).$$

\noind
{\bf (1.2)\quad The $J$-adic filtration and the normal cone.}
Let $J=d({\cal O}^{-1}_X) \subset{\cal O}^0_X={\cal O}_{X^0}$ be the ideal
of the subscheme $\pi_0 (X)$. We denote by
$$N=N_{{\pi_0}(X)/X^0} =\Spec\, \bigoplus\limits_n J^n/J^{n+1}$$
the normal cone of $\pi_0(X)$ in $X^0$. Let also
$$K=\Ker \{d^1:{\bf t}^1_X\to {\bf t}^2_X\}.$$
This is a coherent sheaf on $\pi_0(X)$. Since it
is defined as the kernel of a morphism of
vector bundles, we can associate to it its
total space $\uK\subset \underline{\bf t}^1$ (which is a cone).

\proclaim (1.2.1) Proposition.   There is a natural closed
embedding $N\subset \uK$
of cones over $\pi_0 (X)$.

\prf We have
$$\omega^{-1}=\IO^{-1}={\cal O}^{-1}_X\big|_{{\pi_0}(X)}=
{\cal O}^{-1}_X\big/(d{\cal O}^{-1}_X)\cdot {\cal O}^{-1}_X.$$

\noind
Therefore
$$\underline{\bf t}^1=\Spec\ S(\IO^{-1}),\ \uK=\Spec\ \big(S(\IO^{-1})/
(d\IO^{-2})\cdot S(\IO^{-1})\big)$$
The differential $d^{-1}:{\cal O}^{-1}_X\tworight J$ induces, after
passing to the $n$th symmetric power and restricting
to $\pi_0(X)$, a surjective map
$$\delta_n:S^n(\IO^{-1})\to J^n/J^{n+1}.$$
Explicitly, let $\varphi_1,\cdots, \varphi_n$ be local sections of 
$\IO^{-1}$.
Then
$$\delta_n(\varphi_1\cdots \varphi_n)=d^{-1}(\tilde\varphi_1)\cdots
d^{-1}(\tilde\varphi_n)\,\, {\rm mod}\ J^{n+1}$$
where $\tilde\varphi_i$ is a local section of ${\cal O}^{-1}_X$
extending $\varphi_i$. Therefore
 we get a surjective homomorphism of sheaves
of algebras
$$\delta = \bigoplus \delta_n : S(\IO^{-1})\to \bigoplus J^n/J^{n+1}$$
which induces a closed embedding $\delta^*:N\subset \underline {\bf t}^1$.
To show that $\Im(\delta^*)\subset \uK$, it is enough to
show that $\delta_n((d\varphi)\cdot \varphi_1\cdots\varphi_{n-1})=0$
for any 
local sections $\varphi\in\IO^{-2},\ \varphi_i\in\IO^{-1}$. 
Let $\tilde\varphi\in{\cal O}^{-2}_X,\ \tilde\varphi_i\in{\cal O}^{-1}_X$
be local sections that extend $\varphi, \varphi_i$.
Then,
$$\delta_n\big((d\varphi)\cdot \varphi_1\cdots\varphi_n\big) =
d^{-1}(d^{-2}\tilde\varphi)\cdot d^{-1}(\tilde\varphi_1)\cdots
d^{-1}(\tilde\varphi_n)\,\, {\rm mod}\ J^{n+1}$$
which is clearly $0$.

\noind
{\bf (1.3)\quad Deformation to the normal cone.}
Let $V$ be any vector bundle on $X^0$. We equip it with the $J$-adic filtration
by setting $J^nV = J^n \cdot V$, so that $V$ becomes a filtered module 
over the filtered algebra $({\cal O}_{X^0}, \{J^n\})$. 
Hence ${\rm gr}_J V$ is a graded module over the graded algebra ${\rm gr}_J {\cal O}_{X^0}$
and gives, by localization, a coherent sheaf $\widetilde{\rm gr}_J V$ on
$N = {\rm Spec} \, {\rm gr}_J{\cal O}_{X^0}$. Let $p: N\to \pi_0(X)$
be the projection. The following is well known, 
 with proof supplied for completeness.  

\proclaim (1.3.1) Proposition. The sheaf   $\widetilde{\rm gr}_J V$ is 
identified with $p^*(V\big|_{\pi_0(X)})$. If $f: V\to W$ is any morphism of
vector bundles on $X^0$, then $\widetilde{\rm gr}_J(f)$ is identified with
$p^*(f|_{\pi_0(X)})$. 

\prf Denote for short $Z=\pi_0(X)$. The surjective homomorphism
$V\otimes_{{\cal O}_{X^0}}J^n \to J^n V$ induces, after restricting to $Z$, 
a surjective homomorphism 
$$h_n:(V/JV )\otimes_{{\cal O}_Z} (J^n/J^{n+1}) \to J^nV/J^{n+1}V.$$
We claim that it is an isomorphism. Indeed, if $V= {\cal O}_{X^0}$, the
statement is tautological. Hence it is true for a trivial bundle
$V= {\cal O}^r_{X^0}$. In general, the fact that $h_n$ is an isomorphism
can be verified locally on the Zariski topology, so it follows from local 
triviality of $V$. Now, notice that $V/JV =  V\big|_Z$ and tensoring
with $\bigoplus J^n/J^{n+1}$ over ${\cal O}_Z$ is geometrically
the pullback $p^*$, so $\bigoplus h_n$ gives the required identification.
The statement about morphisms follows from the naturality of maps $h_n$.   

\vskip .2cm

Next, we extend the $J$-adic filtration to ${\cal O}^\b_X$ by setting
$J^n{\cal O}^\b_X$ to have components $J^n\cdot{\cal O}^i_X,\ i\le 0$.
Then $J^n$ is a multiplicative filtration on
the sheaf of dg-algebras ${\cal O}^\b_X$, so ${\rm gr}^\b_J\,{\cal O}^\b_X$
is a graded sheaf of dg-algebras on $X^0$ supported
on $\pi_0(X)$. We have therefore a dg-scheme
$$\Spec\, ({\rm gr}^\b_J\, {\cal O}^\b_X)\lright\pi_0(X).$$\
The underlying ordinary scheme of this dg-scheme
is $\Spec\, ({\rm gr}^\bullet_J\, {\cal O}^0_X)=N$.

Further, let $E^\bullet$ be a dg-vector bundle
on $X$. Then we have the $J$-adic filtration
$J^nE^\bullet$ similarly to the above. The associated
graded object ${\rm gr}^\b_JE^\b$ is then a sheaf
of dg-modules over ${\rm gr}^\b_J\, {\cal O}^\b_X$ and as
such localizes to a sheaf of dg-modules
$\widetilde{\rm gr}_JE^\b$ on the dg-scheme 
$\Spec\, ({\rm gr}_J\, {\cal O}^\b_X)$.
The following is an immediate consequence of Proposition 1.3.1. 

\proclaim (1.3.2) Proposition. (1)
The structure sheaf 
of the dg-scheme $\Spec\, ({\rm gr}^\bullet_J\, {\cal O}^\b_X)$ is
isomorphic, as a sheaf of dg-algebras, to $p^*\IO^\b_X$, where $p^*$
means the
usual pullback of coherent sheaves on 
schemes. \hfill\break
(2) With respect to the identification of (1), 
the sheaf of dg-modules $\widetilde{\rm gr}_JE^\b$ 
is isomorphic to $p^*(E^\b\big|_{{\pi_0}(X)})$.

\proclaim (1.3.3) Proposition. (1) The pullback to $p^*(E^\b|_{\pi_0(X)})$ of the
filtration $\cal D$ is compatible with the differential. \hfill\break
(2) The sheaf of dg-algebras ${\rm gr}_{p^*{\cal D}}(p^*\IO^\b, 
p^*d_{\IO})$ is isomorphic
to $p^*S(\omega^{\le-1})$, the restriction of the
Koszul complex $q^*S(\omega^{\le-1})$ to $N\subset \uK$ (here $q:\uK\rightarrow \pi_0(X)$
is the projection). \hfill\break
(3) The sheaf of dg-modules
 ${\rm gr}_{p^*{\cal D}}\bigl(p^*(E^\b |_{\pi_0(X)})\bigr)$ is isomorphic to
$p^*(\overline{E}\otimes S(\omega^{\leq -1}))$. 

\prf  (1) It is enough to prove that the differential in 
$E^\b|_{\pi_0(X)}$ (denote it $\delta$) 
is compatible
 with the  filtration $\cal D$,  i.e.,
$$\delta \bigl( (\IO^{\leq -1})^n\cdot E^\b |_{\pi_0(X)}\bigr)\i
 (\IO^{\leq -1})^n\cdot E^\b |_{\pi_0(X)}.$$
This follows from the Leibniz rule
$$\delta(fe) = d_{\IO}f\cdot e + (-1)^{{\rm deg}(f) \cdot {\rm deg}(e)}
f\cdot \delta (e), \quad f\in \IO^\b, e\in E^\b |_{\pi_0(X)}$$
and the fact that $d_{\IO}(\IO^{-1})=0$. 
Parts (2) and (3)  follow from (1) and Proposition 1.1.1.

\vskip 2cm

\centerline{\bf 2.\quad Bounded dg-manifolds and [0, 1]-manifolds.}

\vskip 1cm

\noind
{\bf (2.1)\quad $[0, n]$-manifolds.}
Let $X$ be a dg-manifold, and $n\ge 0$.

\proclaim (2.1.1) Proposition. 
The following are equivalent:\hfill\break
{(i)\ }
$\forall x\in\pi_0(X)(\IC)$ the tangent dg-space $T^\b_x X$
is exact outside the degrees in $[0, n]$. \hfill\break
{(ii)\ }
The complex ${\bf t}^\b_X$ is exact outside the
degrees in $[0, n]$.

\prf
(i) $\Rarrow$ (ii) A fiberwise exact sequence of
vector bundles is exact at the level of sheaves
of sections.

(ii) $\Rarrow$ (i) Follows from the spectral sequence
$${\rm Tor}_i^{{\cal O}_{\pi_0(X)}}(\uH^j({\bf t}^\b_X),\
 {\IC}_x)\Rarrow H^{j-i}
(T^\b_x X)$$
and the fact that $T^\bullet_x X$ is situated in
degrees $\ge 0$.

\proclaim (2.1.2)  Definition. We say that $X$ is a  $[0, n]$-manifold
if the conditions of Proposition 2.1.1 are satisfied.

\noindent {\bf (2.1.3)  Examples.}
 (a) If $Y$ is a projective variety of
dimension $n$, then the dg-manifold
${\rm RQuot}_h(\cF)$ constructed in [CK1], 
 is a $[0, n]$-manifold for any
coherent sheaf $\cF$ on $Y$ and any polynomial $h$.

\noindent (b) The dg-manifold ${\rm RHilb}_h^{\rm LCI} (Y)$
 constructed in [CK2],  is a
$[0, d]$-manifold, where $d=\deg(h)$.

\noindent (c) 
Let $X {\buildrel f\over \rightarrow} Z 
{\buildrel g\over\leftarrow} Y$ be a diagram of
smooth algebraic varieties (trivial dg-structure).
Then the derived fiber product 
$X\times_Z^R Y$, constructed in [CK1],  
is a $[0, 1]$-manifold. Indeed, let $(x, y)$ be a
point of
$$\pi_0\biggl(X\times_Z^R Y\biggr) = 
X\times_Z  Y=
\{(x, y)\in X\times Y |f(x)=g(y)\}.$$
and $z=f(x)=g(y)$. Then $T^\b_{(x, y)} 
\bigl(X \times_Z^R Y\bigr)
$ is, up to quasiisomorphism,
 the derived functor  of the fiber product
in the category of vector spaces evaluated on the diagram
$T_xX {\buildrel d_xf\over\longrightarrow} T_zZ 
{\buildrel d_yg\over\longleftarrow}
 T_yY$. This derived functor
is respresented by the 2-term complex
$$T_xX\oplus T_yY\;
{\buildrel d_x f-d_y g\over\longrightarrow}\; T_z Z.$$
In particular, when $f, g$ are closed embeddings,
the derived fiber product is the derived
intersection $X \cap_Z^R  Y$
which is, 
therefore, a $[0, 1]$-manifold.

\vskip .2cm

\noindent {\bf (2.1.4) Remark.} An affine
 $[0,n]$-manifold is the spectrum of
a perfect resolving algebra in the sense of Behrend [B]. 

\vskip .3cm

\noind{\bf 2.2\quad Boundedness and $[0, 1]$-manifolds.}

\proclaim (2.2.1)  Definition.
A dg-manifold $X$ is called
bounded, if $\uH^i({\cal O}^\b_X)=0$ for $i<<0$.

\proclaim (2.2.2) Theorem.
Any $[0, 1]$-manifold is bounded.

\prf
Let $\mu=\max\limits_{x\in\pi_0(X)} \dim\, H^1(T^\b_x X)$.
We will
prove that $\uH^i({\cal O}^\b_X)=0$ for $i<-\mu$.
Since taking cohomology sheaves commutes with
completion, it is enough to prove that $\forall x\in\pi_0(X)(\IC)$\
the complete local dg-ring
$$\widehat{\cal O}^\b_{X, x}={\cal O}^\b_X\otimes_{{\cal O}_X^0}
\widehat{\cal O}_{X^0, x}$$
is exact in degrees $<-\mu$.

\proclaim (2.2.3) Proposition. 
There is a spectral sequence
$$E_2=S^\b\big(H^\b(T^*_x X)\big)\Rarrow H^\b(\widehat{\cal O}^\b_{X, x})$$

\prf
Let $M\subset\widehat{\cal O}^\b_{X, x}$ be the maximal $dg$-ideal
corresponding to $x$, i.e., $M={\bf m}+\widehat{\cal O}^{<0}_{X, x}$ 
where ${\bf m}\subset \widehat{\cal O}_{X^0, x}$ is the usual maximal
ideal in the completed local ring. Then
$$M^n/M^{n+1}\simeq S^n(T^*_x X)$$
as dg-vector spaces, so
$$H^\b(M^n/M^{n+1})=S^n\big(H^\b(T^*_x X)\big).$$
Our spectral sequence is therefore
associated to the filtration $\{M^n\}$.

\vskip .2cm

To finish the proof of Theorem 2.2.2, note that $S^\b(H^\b(T^*_xX))$
is isomorphic to the tensor product of the symmetric algebra of 
$H^0(T^\b_xX)^*$ (situated in degree 0) and the exterior algebra of
$H^1(T^\b_xX)^*$ with the grading being the negative of the usual grading by
the degree of exterior powers. So it clearly vanishes in degrees $<-\mu$.  

\vskip .2cm

\noindent {\bf (2.2.4) Remark.} 
The converse to Theorem 2.2.2 is not true.
For example, if $E$ is a vector bundle on a
manifold $X^0$, then $\big(X^0, \Lambda^\b(E)\big)$ with 
$\deg\, (E)=-3$ 
is bounded but is not a $[0, 1]$-manifold.

\vskip 2cm

\centerline
{\bf 3.\quad The virtual fundamental class of a [0, 1]-manifold.}

\vskip 1cm

\noindent {\bf (3.1) Reminder on Grothendieck and Chow groups.} 
For any  quasiprojective scheme $Y$ we denote by $K_{\circ}(Y)$ the Grothendieck
group of coherent sheaves on $Y$. For such a sheaf $\cal F$ we denote by 
$[{\cal F}]$ its class in $K_\circ (Y)$. We also denote by $K^\circ(Y)$
the Grothendieck ring of vector bundles. As well known, $K_\circ (Y)$
is a module over $K^\circ(Y)$. 
We denote by  $A_r(Y)$  the Chow group of $r$-dimensional cycles on $Y$.
Let $F_r K_\circ(Y)$ be
 the subgroup generated by $[{\cal F}]$ with $\dim \, {\rm supp}\, {\cal F}
\leq r$. Let
$${\rm cl}_r: F_r K_\circ (Y)\to A_r(Y)\otimes {\bf Q},
\quad [{\cal F}]\mapsto \sum_{Z\i {\rm supp}({\cal F}): {\rm dim}(Z)=r}
{\rm mult}_Z({\cal F})\cdot Z$$
be the class map. See [F], Example 18.3.11.  

Let $i: Z\to Y$ be a regular embedding of codimension $d$ such that
 ${\cal O}_Z$ has a finite locally free resolution by ${\cal O}_Y$-modules.
 We denote 
$$i^*_A: A_r(Y)\to A_{r-d}(Z), \quad i^*_K: K_\circ(Y)\to K_\circ(Z)$$
the Gysin maps on the Chow and Grothendieck groups.  Recall that
$$i^*_K([{\cal F}]) = \sum_i (-1)^i [\underline{\rm Tor}_i^{{\cal O}_Y}(
{\cal F}, {\cal O}_Z)].\leqno (3.1.1)$$
Recall also the following ([F], Example 18.3.15).
\proclaim (3.1.2) Proposition. We have
$$i^*_K(F_rK_\circ(Y))\i F_{r-d} K_\circ(Z)$$
and ${\rm cl}_{r-d} i^*_K = i^*_A {\rm cl}_r$.

\vskip .3cm

\noindent {\bf (3.2) The virtual classes.}

\proclaim (3.2.1) Definition.
Let $X$ be a bounded $dg$-manifold.
Its $K$-theoretic virtual fundamental class is defined to be
$[X]^{\rm vir}_{K} =[\uH^\b({\cal O}^\b_X)]\in K_{\circ}\big(\pi_0 (X)\big)$.

From now on we assume that $X$ is a
$[0, 1]$-manifold and use the notation of $\S 1$.

\proclaim (3.2.2) Proposition.
The sheaf $K$ (defined in (1.2)) is locally free.

\prf
This is a consequence of the following lemma.

\proclaim (3.2.3) Lemma. Let $A$ be a Noetherian local ring
with residue field $k$ and
$$Q^1{\buildrel d_1\over\lright}Q^2{\buildrel d_2\over\lright}Q^3$$
an exact sequence of finitely generated free $A$-modules,
which also remains exact after tensoring with $k$. Then
$M=\Ker(d_1)$ is free.

\prf
A finitely generated $A$-module $M$ is free $\iff 
\Tor_1^A(M, k)=0.$

\noind
In our case, the resolution $M\sim\{Q^1{\buildrel d_1\over\lright}
Q^2{\buildrel d_2\over\lright}\ \Im\ d_2\}$
and the fact that $\Ker(d_2\otimes k)=\Im\ (d_1\otimes k)$
implies that $\Tor_1(M, k)=\Tor_{-1}(\Im\ d_2, k)=0$.

\vskip 0.5cm

Since $X$ is a $[0, 1]$-manifold, the truncation
$$\tau_{\le 1}{\bf t}^\b=\{{\bf t}^0\to K\}$$
is quasiisomorphic to ${\bf t}^\b$.

\proclaim (3.2.4) Proposition.
The dual complex $\{K^*\to \omega^0_X\}$
is a perfect obstruction theory on $\pi_0(X)$
in the sense of {\rm {[BF]}}. 

\prf
The embedding of dg-schemes $\pi_0(X)\hookrightarrow X$
induces the morphism of tangent complexes
$$RT^\b\big(\pi_0 (X)\big)\to T^\b X\otimes_{{\cal O}^\b_X}
{\cal O}_{{\pi_0}(X)}={\bf t}^\b$$
Dualizing and passing to truncations, we get
a morphism of 2-term complexes
$$\{K^*\to\omega^0_X\}=\tau_{\ge-1}\omega^\b_X\to\tau_{\ge-1}
L\Omega^{1\b}\big(\pi_0(X)\big)
\cong\{J/J^2\to\Omega^1_{X^0}\big|_{{\pi_0}(X)}=\omega^0\}$$
which is clearly an isomorphism on $\uH^0$. 

\noind
Explicitly, this morphism is identical on the $0$th
terms and on the $(-1)$st terms is induced by
the surjective map
$$d:{\cal O}^{-1}/d{\cal O}^{-2}\longrightarrow J=d{\cal O}^{-1}$$
after restricting to $\pi_0$. So we have a morphism
of $2$-term complexes which is an isomorphism
in degree $0$ and a surjection in degree $(-1)$,
inducing an isomorphism on $\uH^0$ and surjection
on $\uH^{-1}$. This is precisely the definition of a
perfect obstruction theory.

\vskip 0.5cm

Following [BF], we give

\proclaim (3.2.5) Definition. Let $i: \pi_0(X)\to \underline{K}$ be the
embedding of the zero section. 
The homological virtual fundamental
class of $X$ is the element
$$[X]^{\rm vir}= i^*_A [N] \in A_r\big(\pi_0(X)\big).$$
Here $r={\rm vdim} (X)={\rm rk}({\bf t}^0)-{\rm rk}(K)$ is the
virtual dimension of $X$.

\vskip .3cm

\proclaim (3.3) Theorem.
The $K$-theoretic fundamental class $[X]^{\vir}_K$ lies
in $F_r K\big(\pi_0 (X)\big)$ and 
$${\rm cl}_r\big([X]^\vir_K\big)=[X]^\vir .$$

\prf By (3.1.1-2), 
it is enough to show that
$$[X]^\vir_K=\sum(-1)^i[\underline\Tor\, _i^{{\cal O}_\uK}
({\cal O}_N, {\cal O}_{{\pi_0}(X)})]\in K_{\circ}\big(\pi_0(X)\big)$$
the sum being finite since ${\cal O}_{{\pi_0}(X)}$ has a finite
locally free resolution over ${\cal O}_{\uK}$, namely the
Koszul complex. Denoting $q:\uK\to\pi_0(X)$
the projection, we can write the Koszul
resolution as $\Lambda^\b(q^*K^*)\sim{\cal O}_{{\pi_0}(X)}$, with
the differential induced by the tautological
section $\xi\in \Gamma(\uK, q^*K)$. The embedding
$K\subset {\bf t}^1$ defines a quasiisomorphism
$$\varphi:K\to {\bf t}^{\ge 1}[1]=\{{\bf t}^1\to {\bf t}^2\to\cdots\},\ 
\deg({\bf t}^i)=i-1.$$
In particular, we have a section $q^*(\varphi)(\xi)$ of
the dg-bundle $q^*({\bf t}^{\ge 1}[1])$ on $\uK$ and the
induced Koszul complex $q^*\big(S(\omega^{\le-1}_X)\big)$ is
a resolution of ${\cal O}_{{\pi_0}(X)}$ on $\uK$.

The direct image map $i_*:K_{\circ}\big({\pi_0}(X)\big)\to K_{\circ}(\uK)$
preserves the dimension
filtration.
By the above discussion, if $i_*$ were injective, our theorem would follow from
the equality
$$i_*[\uH^\b({\cal O}^\b_X)]=[\uH^\b(q^* S(\omega^{\le-1})
\otimes_{{\cal O}_\uK}{\cal O}_N)]\leqno (3.3.1)$$
in $K_\circ(\uK)$.
While, in general, $i_*$ is not injective,
it becomes so after passing to {\it equivariant} $K$-theory.
Specifically, consider the $\Bbb G_m$-action on $\uK$ given by dilations on the fibers.
By a slight abuse of notation, let us denote by
$$i_*:K_{\circ}^{\Bbb G_m}\big({\pi_0}(X)\big)\to K_{\circ}^{\Bbb G_m}(\uK)$$
the direct image map in $\Bbb G_m$-equivariant $K$-theory.
Recall that
$$K_\circ^{\Bbb G_m}({\rm pt}) = \Bbb C[\mu, \mu^{-1}].$$
The section $i$ embedds $\pi_0(X)$ into $\uK$ as the fixed point locus of the action.
Therefore,
it follows from the localization theorem ([T], Thm. 2.1)
that $i_*$
becomes an isomorphism after tensoring with the quotient field $\Bbb C(\mu)$
of $\Bbb C[\mu, \mu^{-1}]$.
Since
$$K_{\circ}^{\Bbb G_m}\big({\pi_0}(X)\big)\cong K_{\circ}\big({\pi_0}(X)\big)\otimes \Bbb C[\mu,\mu^{-1}]$$
has no $\Bbb C[\mu,\mu^{-1}]$-torsion, we conclude that
(the equivariant version of) $i_*$ is injective.
Further, if we consider $K, {\bf t}^1, {\bf t}^2,...$
as equivariant bundles on $\pi_0(X)$ (with $\Bbb G_m$ acting by dilations
in the fibers) and use the equivariant flat
pull-back $q^*$, then the Koszul complexes $\Lambda^\b(q^*K^*)$ and
$q^*\big(S(\omega^{\le-1}_X)\big)$
are equivariant resolutions of ${\cal O}_{\pi_0(X)}$.
Finally, we have the $\Bbb G_m$-equivariant Gysin map
$i^*$ (defined by the same formula with Tor's) which satisfies
$$i_*(i^*({\cal F}))=\Lambda^\b(q^*K^*)\otimes_{{\cal O}_\uK}{\cal F},\;\;\;
{\cal F}\in K_{\circ}^{\Bbb G_m}(\uK).$$
We conclude that it is indeed sufficient to prove the equality (3.3.1), but in
an upgraded form, in
which all maps and sheaves are  considered in $\Bbb G_m$-equivariant $K$-theory.
So in the rest of the proof we will deal with equivariant theory.

Let us factor $i$ into the
composition of two embeddings
$$\pi_0(X){\buildrel i_3\over \hookrightarrow}N
{\buildrel i_2\over \hookrightarrow}\uK.\leqno (3.3.2)$$
The inclusions $i_2$, $i_3$,
as well as the projection $p:N\to \pi_0(X)$, are
$\Bbb G_m$-equivariant and we use the corresponding
equivariant push-forward or pull-back maps.

Note that the RHS of (3.3.1) is equal to
$$i_{2*}[\uH^\b(p^* S\big(\omega^{\le 1})\big)].\leqno(3.3.3)$$

Recall (Proposition 1.3.2) that the
sheaf of dg-algebras $p^*\IO^\b$ on $N$ is the localization on
$N=\Spec\ {\rm gr}_J {\cal O}_{X^0}$ of the sheaf of graded
$dg$-algebras ${\rm gr}_J\ {\cal O}^\b_X$.
Proposition 1.3.3(2)  implies that
$$[\uH^\b(p^*\IO, \delta)]=[\uH^\b\big(p^*S(\omega^{\le-1})\big)]
\in K_\circ^{\Bbb G_m}(N)\leqno (3.3.4)$$
by virtue  of the spectral sequence of the filtered complex
$\big((p^*\IO, \delta), p^*{\cal D}\big)$.

Next, we have a spectral sequence of
coherent sheaves on $N$
$$\uH^\b(p^*\IO, \delta)\Rightarrow i_{3*}\uH^\b({\cal O}^\b_X).$$
It is obtained from Proposition 1.3.2 
by localizing the spectral sequence of
the sheaf of filtered $dg$-algebras $({\cal O}^\b_X, J)$
over $N=\Spec\ {\rm gr}_J\ {\cal O}_{X^0}$. This spectral sequence
converges by the Noetherian argument and 
implies the equality
$$[\uH^\b(p^*\IO,\delta)]=i_{3*}[\uH^\b({\cal O}^\b_X)]\in K_\circ^{\Bbb G_m}(N)
\leqno(3.3.5)$$
Putting together (3.3.3), (3.3.4) and (3.3.5), we get (3.3.1).

\vskip 2cm

\centerline
{\bf 4.\quad The virtual class via the Chern character.}

\vskip 1cm

\noind
{\bf 4.1\quad Reminder on local Chern character and Riemann-Roch.}
Let $Z\to Y$ be a closed embedding of schemes
of finite type over $k$. We denote by $A^m(Z\to Y),\
m\in\IZ$, the $m$th operational Chow group [F].
Its elements act by homomorphisms $A_p(Y)\to A_{p-m}(Z)$,
and $A_p(Z)=A^{-p}(Z\to pt)$. When $Y$ is smooth,
$A^m(Z\to Y)=A_{\dim Z-m}(Z).$

If $F^\b$ is a finite complex of vector bundles
on $Y$ exact outside of $Z$, one has the
localized Chern character
$${\rm ch}^Y_Z(F^\b)\in A^\b(Z\to Y)\otimes \IQ.$$

\noind
We denote by
$$\tau_Z:K_\circ(Z)\to A_\b (Z)$$
the Riemann-Roch map of Baum-Fulton-McPherson
[F]. If $Y$ is a smooth quasiprojective
variety containing $Z$ as a closed subscheme,
and $\cF$ is a coherent sheaf on $Z$, then
$$\tau_Z[\cF]={\rm ch}^Y_Z(F^\b)\cdot {\rm Td}(T_Y)\leqno(4.1.1)$$
where $F^\b$ is a locally free resolution
of $\cF$ on $Y$.

For any proper morphism $f:Z\to W$ 
of quasiprojective schemes we denote
$$f^A_*:A_\b(Z)\to A_\b(W)$$
the direct image map on the Chow groups.
The Riemann-Roch theorem in the form
of Baum-Fulton-McPherson (see [F], Th. 18.2) says that
$$\tau_W\big(f_*(z)\big)=f^A_*\big(\tau_Z(z)\big),\
z\in K_\circ(Z).\leqno(4.1.2)$$

Let $Z\to Y$ be a closed embedding of
quasiprojective schemes and $F^\b$ a finite
complex of vector bundles on $Y$, exact outside
$Z$. Then for any coherent sheaf $\cal G$ on $Y$
we have the Riemann-Roch formula ([F], Ex. 18.3.12):
$$\tau_Z [\underline{H}^\b(F^\b\otimes{\cal G}^\b)] = {\rm ch}_Z^Y(F^\b)\cdot 
\tau_Y({\cal G}). \leqno (4.1.3)$$
Let now $Z$ be proper.
Combining (4.1.3) for ${\cal G}={\cal O}_Y$ with the
formula (4.1.2) for $W={\rm pt}$, we get the
following form of the Riemann-Roch
theorem:
$$\chi \big(Y, \uH^\b(F^\b)\big)=\int\limits_Z {\rm ch}^Y_Z(F^\b)\cdot {\rm Td}
(T_Y).
\leqno (4.1.4)$$

\vskip .3cm

\noind
{\bf (4.2)\quad Kontsevich's definition of the homological virtual
class.}
Let $X$ be a $[0, 1]$-manifold. In [K], M. Kontsevich
proposed to consider the element
$$\kappa_X=\tau_{{\pi_0}(X)}[\uH^\b({\cal O}^\b_X)]\cdot {\rm Td}^{-1}
({\bf t}^\b_X)\in A_\b\big(\pi_0(X)\big)\otimes \IQ\leqno (4.2.1)$$ 
as the virtual fundamental class of $X$.
Since we use the embedding $\pi_0(X)\subset X^0$ 
for the definition of $\tau_{{\pi_0}(X)}$, we have, applying
(4.1.3) to $F^\b={\cal O}^\b_X,\ {\cal G} ={\cal O}_{X^0}$, that
$$\kappa_X={\rm ch}^{X^0}_{{\pi_0}(X)}({\cal O}^\b_X)\cdot
 {\rm Td}({\bf t}^{\ge 1}_X
[1])\leqno (4.2.2)$$

\proclaim (4.2.3) Theorem.
$\kappa_X=[X]^\vir$ (equality in $A_\b(\pi_0(X))\otimes \IQ$).
 In particular, $\kappa_X$ is
homogeneous of degree ${\rm vdim}(X)$.

\noind
{\bf (4.3)\quad Proof of Theorem 4.2.3.}
We use
the notation introduced in the proof of Theorem 3.3,
in particular, the embeddings $i_2, i_3$, their composition $i$ and the
projections $p, q$. We need to show that $\kappa_X=i^*_A([N])$. 

Using the quasiisomorphism
${\bf t}^\b_X\sim\{{\bf t}^0_X\to K\}$, we have
$$\kappa_X=\tau_{\pi_0(X)}\big([\uH^\b({\cal O}^\b_X)]\big)\cdot {\rm Td}^{-1}
({\bf t}^0_X)\cdot {\rm Td}(K).$$

We have shown in Theorem 3.3 that $[\uH^\b({\cal O}^\b_X)]=i^*({\cal O}_N)$. On the other hand,
since $i$ is a regular embedding with normal bundle $K$, we have by [F, Thm. 18.2 (3)]
$$\tau_{\pi_0(X)}\big(i^*({\cal O}_N)\big)={\rm Td}^{-1}(K)\cdot i^*_A\big(\tau_{\uK}({\cal O}_N)\big),$$
hence
$$\kappa_X=i^*_A\big(\tau_{\uK}({\cal O}_N)\big)\cdot {\rm Td}^{-1}
({\bf t}^0_X).$$ 
Now use the Riemann-Roch formula (4.1.2) for the embedding $i_2$ to get
$$\kappa_X=i^*_A\big(i_{2*}^A\tau_{N}({\cal O}_N)\big)\cdot {\rm Td}^{-1}
({\bf t}^0_X).\leqno(4.3.1)$$ 
Our proof is then a consequence of the following 

\proclaim (4.3.2) Lemma.
Let $i:Z\subset Y$ be a closed embedding
of quasiprojective schemes with $Y$ smooth
and $p:N\to Z$ be the projection of the
normal cone $N=N_{Z/Y}$. Then
$$\tau_N({\cal O}_N)=p^* {\rm Td}(T_Y\big|_Z)\cdot [N].$$

Specifically, we apply the lemma to
$Z=\pi_0(X),\
Y=X^0$, so that ${\bf t}^0_X=T_Y\big|_Z$, getting
from (4.3.1) that
$$\kappa_X=i^*_A\big(i_{2*}^A(p^*{\rm Td}({\bf t}^0_X)\cdot [N])\big)\cdot {\rm Td}^{-1}
({\bf t}^0_X).$$
Since $p^*=i_2^*q^*$ (with $i_2^*$ the pull-back on operational Chow rings), the projection formula gives
$$\kappa_X=i^*_A\big(q^*{\rm Td}({\bf t}^0_X)\cdot i^A_{2*}([N])\big)\cdot {\rm Td}^{-1}
({\bf t}^0_X).$$ But the right-hand side of the last equality is 
precisely $i^*_A([N])$, as $q\circ i=id_{\pi_0(X)}$.

\vskip .3cm

\noind
{\bf Proof of Lemma 4.3.2.} Let $J\subset {\cal O}_Y$ be the ideal
of $Z$ and
$$\tilde Y=\Spec\bigoplus^\infty\limits_{n=0} J^n\cdot t^n
{\buildrel \varrho\over\lright} \IA^1=\Spec\ \IC[t]$$
be the deformation to the normal cone. The morphism
$\varrho$ is flat, with $\varrho^{-1}(0)=N$ and
$\varrho^{-1}(t)\simeq Y,\ t\not=0$.
Let $\varepsilon_t:\varrho^{-1}(t)\hookr\tY$ be the embedding
and 
$$\varepsilon_t^!:A_\b(\tY)\to A_{\b -1}(\varrho^{-1}(t))$$
be the specialization map of [F],  10.1. By
Example 18.3.8 of [F]
$$\tau_{\varrho^{-1}(t)} ({\cal O}_{\varrho^{-1}(t)})=
\varepsilon^{!}_t \tau_{\tY}({\cal O}_{\tY}),\ t\in \IA^1.$$
Moreover, $\tau_{\tY}({\cal O}_{\tY})$ is uniquely defined by its
specializations for $t\not=0$ ([F], 11.1).
In other words, if $y\in A_\b (\tY)$ is such that
$\varepsilon^!_t(y)=\tau_Y({\cal O}_Y),\ t\not=0$, then necessarily
$y=\tau_{\tY}({\cal O}_{\tY})$
and hence $\tau_N({\cal O}_N)=\varepsilon^{!}_0(y)$.

We have a projection $\sigma:\tY\to Y$ 
induced by the embedding ${\cal O}_Y=J^0\cdot t^0\subset
\bigoplus J^n\cdot t^n$.
The map $\sigma$ is the identity on each $\varrho^{-1}(t)=Y,\ t\not=0$
and is equal to $ip$ on $\varrho^{-1}(0)=N$. Let now
$y=\big(\sigma^* {\rm Td}(T_Y)\big)[\tY]\in A_\b(\tY)$. Here we view
${\rm Td}(T_Y)$ as an element of $A^\b(Y)$, so $\sigma^*{\rm Td}(T_Y)\in
A^\b(\tY)=A^\b(\tY\to \tY)$ and $y$ is the value of
$\sigma^*{\rm Td}(T_Y)$ on $[\tY]\in A_\b (\tY)$. Then, clearly, $y$ satisfies
the above condition on $\varepsilon^!_t(y),\ t\not=0$, so
$$\tau_N({\cal O}_N)=\vep^!_0(y)=\vep^!_0\big(\sigma^*{\rm Td}(T_Y)\big)[\tY]=
p^*i^*{\rm Td}(T_Y)[N]$$
as claimed.
\vskip .3cm

\noind
{\bf (4.4)\quad A Riemann-Roch theorem for dg-manifolds.}
Let $X$ be a $[0, 1]$-manifold. A dg-vector
bundle $E^\b$ on $X$ will be called finitely generated,
if the complex $\oE^\b$ of vector
bundles on $\pi_0(X)$, see (1.1), is finite.
In this case $\uH^j(E^\b)=0$ except for finitely
many $j$ and so we have the class
$[\uH^\b(E^\b)]\in K_\circ\big(\pi_0(X)\big)$.

\proclaim (4.4.1) Theorem.
$$\tau_{{\pi_0}(X)}[\uH^\b(E^\b)]={\rm ch} (\overline{E}^\b)\cdot {\rm Td}
({\bf t}^\b_X)\cdot[X]^\vir.$$
Here the first two factors on the right
are considered as endomorphisms of $A_\b\big({\pi_0}(X)\big)\otimes \IQ$
and applied successively to $[X]^\vir$.

This is a consequence of (4.1.3), of Theorem 4.2.3, and the folllwing fact.

\proclaim (4.4.2) Theorem. We have the equality in $K_\circ(\pi_0(X))$:
$$[\underline{H}^\bullet(E^\b)] = [\overline{E}^\b]\cdot [\uH^\bullet(\cO_X^\b)].$$
(Product of an element of $K^\circ$ with an element of $K_\circ$.)

\prf
We use the equivariant set-up and the notation from the proof of 
Theorem 3.3. Since the $\Bbb G_m$-equivariant push-forward $i_*=i_{2*}i_{3*}$
is injective, it is enough to show that
$$i_{3*}[\uH^\b(E^\b)] = i_{3*}\bigl( [\overline{E}^\b] \cdot [\uH^\b(\cO_X^\b)]\bigr). 
\leqno (4.4.3)$$
This would follow
if we proved the following equality in
$K_\circ^{\Bbb G_m}(N)$:
$$i_{3*}\big[\uH^\b(E^\b)]=[p^*\oE^\b\otimes \Lambda^\b(p^*K^*)\big].
\leqno (4.4.4)$$
The proof of (4.4.4) proceeds similarly to the
case $E^\b={\cal O}^\b_X$, see (3.3.4-5).
To be precise, $\Lambda^\b(p^*K^*)$ has the Koszul differential,
so the RHS of (4.4.4) is equal to
$$\Big[\uH^\b\big(p^*\oE^\b\otimes\Lambda^\b(p^*K^*)\big)\Big]$$
which, in view of the quasiisomorphism $K\to {\bf t}^{\ge 1}[1]$
gives
$$[p^*\oE^\b\otimes \Lambda^\b(p^*K^*)\big]=
\Big[\uH^\b(p^*(\oE^\b\otimes S(\omega^{\le-1})))\Big]. \leqno (4.4.5)$$
By Proposition 1.3.2(2),
$$p^*(E^\b\big|_{{\pi_0}(X)})\simeq \widetilde{\rm gr}_J\ E^\b.\leqno (4.4.6)$$
On the other hand, by Proposition 1.3.3,
$${\rm gr}_{p^*\cD} p^*(E^\b\big|_{{\pi_0}(X)})\simeq p^*\oE^\b\otimes
p^*S(\omega^{\le-1}). \leqno (4.4.7)$$
The spectral sequence of the filtered complex
 $(p^*(E^\b|_{\pi_0(X)}), p^*{\cal D})$ (together with finite generation of 
$E^\b$) implies then that
$[\uH^\b p^*(E^\b\big|_{{\pi_0}(X)})]$ makes sense
and
$$[\uH^\b p^*(E^\b\big|_{{\pi_0}(X)})]=
[\uH^\b(p^*\oE^\b\otimes p^*S(\omega^{\le-1}))].\leqno (4.4.8)$$
Next, (4.4.7) and the spectral sequence of the filtered complex
$(E^\b, J)$ implies 
$$[\uH^\b(p^*E^\b\big|_{{\pi_0}(X)})]=i_{3*}[\uH^\b(E^\b)].\leqno (4.4.9)$$
Combining (4.4.5), (4.4.8) and (4.4.9) proves the
equality (4.4.4) and therefore Theorems 4.4.2 and  4.4.1.

\vskip .2cm

\proclaim (4.4.10) Corollary. For two finitely generated dg-bundles $E^\b, F^\b$ on $X$
we have the equality in $K_\circ(\pi_0(X))$:
$$\bigl[ \uH^\b\bigl( E^\b\otimes_{\cO_X^\b} F^\b\bigr)\bigr] = 
[\overline{E}^\b]\cdot [\uH^\b(F^\b)].$$

\vskip .3cm

\noind
{\bf (4.5)\quad Consequences for the Euler characteristic.}
Let us assume, in the situation of (4.4) that
$\pi_0(X)$ is projective. Then the
Euler characteristic
$$\chi\big(\pi_0(X),\ \uH^\b(E^\b)\big)=\sum(-1)^i\chi\big(\pi_0(X),\
\uH^i(E^\b)\big)$$
is defined. Theorem 4.4.1 allows us to establish 
a simple formula for this Euler characteristic.

Since $\pi_0(X)$ is projective, we have the degree
map
$$\deg:A_0\big(\pi_0(X)\big)\to \IZ.$$
For any $\varphi\in A^\b\big(\pi_0(X)\big)=
A^\b\big(\pi_0(X)\to \pi_0(X)\big)$
we define
$$\int\limits_{[X]^\vir} \varphi=\deg\big((\varphi\cdot[X]^\vir)_0\big) .$$
Here the subscript $0$ means the degree $0$ component
of $\varphi\cdot[X]^\vir\in A_\b\big(\pi_0(X)\big)$.

\proclaim (4.5.1) Theorem.
$$\chi\big({\pi_0}(X),\ \uH^\b(E^\b)\big)=\int\limits_{[X]^\vir}
{\rm ch} (\oE^\b)\cdot {\rm Td}({\bf t}^\b_X).$$

\prf
This is a direct consequence of
Theorem 4.4.1 and the fact that $\tau$ commutes with
direct image (for the map $\pi_0(X)\to {\rm pt}$).

\vskip .3cm

\noindent {\bf (4.6) Chern numbers and the cobordism class
of a [0,1]-manifold.} Let $X$ be a $[0,1]$-manifold of virtual dimension $d$.
Let $P(d)$ be the set of partitions of $d$ into ordered summands, i.e.,
of sequences $I=(i_1, ..., i_p)$ with $i_\nu\in \Bbb Z_+$ and $\sum i_\nu = d$.
For each $I\in P(d)$ we define the $I$th Chern number of $X$ to be
$$c_I(X) = \int_{[X]^{\vir}} c_{i_1}({\bf t}^\b_X) ... c_{i_p}({\bf t}^\b_X) \quad\in\quad \Bbb Z.
\leqno (4.6.1)$$
Let $\Omega U^d$ be the cobordism group of compact almost complex manifolds of
real dimension $2d$, see [R]. For each such manifold $M$ the tangent bundle $T_M$
is a complex vector bundle so it has Chern classes $c_i(T_M)\in H^{2i}(M, \Bbb Z)$,
and for each $I\in P(d)$ we have the Chern number
$$c_I(M) = \int_{[M]} c_{i_1}(T_M) ... c_{i_p}(T_M) \quad\in\quad \Bbb Z. \leqno (4.6.2)$$
Here $[M]$ is the usual fundamental class of $M$. The following is well known, see [R]:

\proclaim (4.6.3) Proposition. (a) The Chern numbers are cobordism invariant.\hfill\break
(b) If two almost complex manifolds have the same Chern numbers, then they are cobordant.

Our next result shows that a $[0,1]$-manifold over $\Bbb C$ can be seen as a ``virtual''
smooth complex manifold. This agrees with the intuition that working with dg-manifolds
is a replacement of deforming to transverse intersection, a technique that typically
leads outside of algebraic geometry.

\proclaim (4.6.4) Theorem. Let $k=\Bbb C$ and $X$ be a $[0,1]$-manifold over $\Bbb C$ of
virtual dimension $d$. Then there exists a (unique, up to cobordism) almost complex
manifold $M$ of real dimension $2d$ such that $c_I(M)=c_I(X)$ for all $I\in P(d)$.

\prf We first recall the concept of Schur functors [M]. Let $\alpha = (\alpha_1 \geq \alpha_2 \geq ...)$
be a weakly decreasing sequence of nonnegative integers terminating in zeroes.
Let also ${\rm Vect}_k$ be the category of finite-dimensional $k$-vector spaces. 
Then we have the Schur functor $\Sigma^\alpha: {\rm Vect}_k\to {\rm Vect}_k$.
If $V=k^d$, then $\Sigma^\alpha(k^d)$ is ``the'' space of the irreducible 
representation of the algebraic group $GL_d/k$ with highest weight $\alpha$.
The functor $\Sigma^\alpha$ can be applied to vector bundles (and projective
modules over any commutative $k$-algebra). In particular, if $k=\Bbb C$
and $M$ is an almost complex manifold, then we have the complex
vector bundle $\Sigma^\alpha(T_M)$ on $M$. In this case the number
$$\phi_\alpha(M) = \int_{[M]} {\rm ch}(\Sigma^\alpha(T_M))\cdot {\rm Td}(T_M)$$
is expressible as a universal $\Bbb Q$-linear combination of the Chern numbers of $M$:
$$\phi_\alpha(M) = \sum_I q_\alpha^I c_I(M), \quad q_\alpha^I\in\Bbb Q. 
\leqno (4.6.5)$$
The following is a reformulation of the Hattori-Stong theorem, see [R] [S]:

\proclaim (4.6.6) Theorem. Let $(\lambda_I)_{I\in P(d)}$ be a system of integers
labelled by $P(d)$. Then the following are equivalent:\hfill\break
(i) There exists an almost complex manifold $M$ (unique up to cobordism)
such that $c_I(M)=\lambda_I$ for all $I\in P(d)$. \hfill\break
(ii) For any $\alpha$ as above the number $\sum_I q^I_\alpha \lambda_I$ is an integer. 

We now prove that the condition (ii) holds for $\lambda_I=c_I(X)$ where $X$ is 
a $[0,1]$-manifold of virtual dimension $d$. Indeed, the Schur functors apply equally
well to dg-bundles on $X$. See, e.g.,   [ABW] for Schur functors of complexes.
If $E^\b$ is a finitely generated bundle, then so is $\Sigma^\alpha(E^\b)$.
Further, Schur functors commute with restrictions of bundles, so in particular,
$$\overline{\Sigma^\alpha(E^\b)} = \Sigma^\alpha(\overline{E}^\b).$$
Now, applying Theorem 4.5.1, we see that
$$\sum_I q^I_\alpha c_I(X)\quad  = \quad \int_{[X]^{\vir}} {\rm ch}(\Sigma^\alpha {\bf t}^\b_X) 
\cdot {\rm Td}({\bf t}^\b_X) \quad = \quad \chi (X, \Sigma^\alpha T^\b_X)\quad\in\quad \Bbb Z,$$
whence the statement.

\vskip 2cm

\centerline{\bf 5.\quad Localization}

\vskip 1cm

\noind
{\bf (5.1) Background.} 
Let $G=(\Bbb G_m)^n$ be an $n$-dimensional
algebraic torus over $k$. For a $G$-scheme $Z$ we
denote by $K_G(Z)$ the Grothendieck group of
$G$-equivariant coherent sheaves on $Z$ and by
$K_G^{\circ}(Z)$ the Grothendieck ring of $G$-vector bundles
on $Z$. We denote by ${\rm Rep}(G)=K_G(pt)$ the
representation ring of $G$ (which is a Laurent 
polynomial ring) and by ${\rm FRep}(G)$ its field of fractions.

\proclaim (5.1.1) Lemma.
If the $G$-action on $Z$ is trivial, and $Z$ is quasiprojective, then,
for every $G$-bundle $E$ satisfying $E^G=0$, the element $[\Lambda^\b(E)]$
is invertible in the localization $K^{\circ}_G(Z)\otimes_{{\rm Rep(G)}} 
{\rm FRep}(G)$.

Let $Y$ be a smooth quasiprojective $G$-variety
and $Z\subset Y$ an invariant closed subscheme.
We will
need a version of the Bott localization formula
for $Z$.

 Denote $\epsilon:Z^G\to Z,\ \widetilde {\epsilon}:Y^G\to Y$ the
embeddings of the fixed point loci, so we have the Cartesian square
of closed embeddings:
$$\matrix{&Z&\buildrel j\over\longrightarrow &Y&\cr
\epsilon&\uparrow&&\uparrow&\widetilde{\epsilon}\cr
&Z^G&\buildrel \widetilde{j}\over\longrightarrow&Y^G&}\leqno (5.1.2)$$
 Note
that $\widetilde{\epsilon}$ is a regular embedding (and  $Y^G$ is
smooth).  Let $\cal N$ be the normal
bundle of $Y^G$ in $Y$ and ${\cal N}^*$ its dual bundle. Let
$$\epsilon^!: K_G(Z)\to K_G(Z^G)\leqno (5.1.3)$$
be the K-theoretic Gysin map defined by putting, for each coherent $G$-sheaf
$\cal F$ on $Z$:
$$\epsilon^!([{\cal F}]) = \sum_l (-1)^l \bigl[ \widetilde{j}^* \underline{\Tor}_l^{\cO_Y}
(j_*{\cal F}, \cO_{Y^G})\bigr]. \leqno (5.1.4)$$
Here the Tor-sheaves are supported on $Z^G$. This is a K-theoretic analog of
the refined Gysin map of Fulton [F]. Like that map, $\epsilon^!$ depends
not only on the morphism $\epsilon$, but on the entire diagram (5.1.2). 

\proclaim (5.1.5) Theorem. For any $\xi\in K_G(Z)$ we have the equality
 $$\xi =\epsilon_*\left({\epsilon^! (\xi)\over
\left[\Lambda^\b\big({\cal N}^*\big|_{Z^G}\big)\right]}
\right)$$
in the group
$K_G(Z)\otimes_{{\rm Rep} (G)} {\rm FRep}(G)$. 

\prf By the result of Thomason ([T], Th. 2.1),
$$\epsilon_*: K_G(Z^G)\otimes_{{\rm Rep} (G)} {\rm FRep}(G)\to K_G(Z)\otimes_{{\rm Rep} (G)} {\rm FRep}(G)
\leqno (5.1.6)$$
is an isomorphism, so $\xi=\epsilon_*(\eta)$  for some $\eta$ in the LHS of (5.1.6). 
On the other hand, for any cohereht $G$-sheaf $\cL$ on $Z^G$ we have
$$\epsilon^! \epsilon_*[\cL] = 
\bigl[ \underline{\Tor}_\b^{\cO_Y}(\widetilde{j}_*\cL, \cO_{Y^G})\bigr]=
\bigl[ \widetilde{j}_*\cL\otimes_{\cO_{Y^G}}\underline{\Tor}_\b^{\cO_Y}(\cO_{Y^G}, \cO_{Y^G})\bigr]
=[\cL]\cdot \bigl[ \Lambda^\b\bigl({\cal N}^*\big|_{Z^G}\bigr)\bigr].$$
Therefore
$$\epsilon^!\xi = \eta\cdot \bigl[ \Lambda^\b\bigl({\cal N}^*\big|_{Z^G}\bigr)\bigr].$$
This means that the fraction in the RHS of the equality claimed in Theorem 5.1.5,
is equal to $\eta$, and the equality is true since $\xi = \epsilon_*(\eta)$.

\vskip .3cm

\noindent{\bf (5.2)\quad The setup.}
Let $X$ be a $[0, 1]$-manifold 
with $G$-action. Then we have the fixed point
(dg-)submanifold $X^G\subset X$, with
$$\eqalign{
(X^G)^0&=(X^0)^G,\ \pi_0(X^G)=\pi_0(X)^G,\cr
{\cal O}^\b_{X^G}&=\big({\cal O}^\b_X\big|_{(X^0)^G}\big)_G \;\; {\rm
(the\; coinvariants)}.\cr}$$

\noind
Let $i:X^G\hookr X$ be the embedding and
$\nu^\b=i^*T^\b_X/ T^\b_{X^G}$ be the dg-normal bundle
of $X^G$. It has the induced $G$-action.
As in (1.1) we denote by $\overline\nu\,^\b$ the restriction
of $\nu^\b$ to $\pi_0(X)^G$ in the sense of dg-manifolds.
Thus we have a split exact sequence of complexes
of vector bundles
$$0\to {\bf t}^\b_{(X)^G}\to {\bf t}^\b_X\big|_{{\pi_0}(X)^G}
\to \overline\nu\,^\b\to 0,\;\;\;\;
{\bf t}^\b_{X^G}=\Big({\bf t}^\b_X\big|_{{\pi_0}(X)^G}\Big)^G\leqno (5.2.1)$$
It shows the following

\proclaim (5.2.2.) Proposition.
$X^G$ is again a $[0, 1]$-manifold, and
$\overline\nu\,^\b$ is fiberwise exact outside of degrees $0, 1$.

Therefore
$$\overline\nu\,^{\b}{}'=\Big\{\overline\nu\,^0\to\Ker\{\overline\nu\,^1
{\buildrel d\over\lright \overline\nu\,^2}\}\Big\}\leqno (5.2.3)$$
is a $2$-term $G$-complex of bundles on $\pi_0(X)^G$
quasiisomorphic to $\overline\nu\,^\b$. This is precisely the
``moving part'' of the obstruction theory ${\bf t}^\b_X$ 
in the sense of [GP].

\vskip .3cm

\noind
{\bf (5.3)\quad $K$-theoretic localization for $[0, 1]$-manifolds.}
In the setup of (5.2) let $E^\b$ be a finitely
generated $G$-equivariant dg-vector bundle on $X$.
We denote by $i^*E^\b=(i^0)^{-1}E^\b\otimes_{(i^0)^{-1}{\cal O}^\b_X}
{\cal O}^\b_{X^G}$
the restriction of $E^\b$ to $X^G$ in the sense of
dg-manifolds. We have the class $[\uH^\b(E^\b)]\in K_G\big(\pi_0(X)\big)$.
In particular, for $E^\b={\cal O}^\b_X$ we get
the $G$-equivariant version of the $K$-theoretic
virtual class
$$[X]^{\vir, G}_K =[\uH^\b({\cal O}^\b_X)]\in K_G\big(\pi_0(X)\big),$$
and, furthermore,
$$[\uH^\b(i^*{\cal O}^\b_X)]=[X^G]^{\rm vir}_K.$$

\proclaim (5.3.1) Theorem.
In $K_G\big(\pi_0(X)\big)\otimes_{{\rm Rep}(G)} {\rm FRep}(G)$ we have
the equality
$$[\uH^\b(E^\b)]=\pi_0(i)_*\left({\Big[\uH^\b(i^*E^\b)\Big]\over
\Big[\Lambda{}{{} ^{\b}}(\overline\nu^\b{}'{}^*)\Big]}\right), $$
where $[\Lambda^\b(\overline\nu^\b{}'{}^*)]$ is defined as 
$[\Lambda^\b(\overline\nu^0{}'{}^*)]/[\Lambda^\b(\overline\nu^1{}'{}^*)]$,
see { (5.2.3)}.

\prf We apply Theorem 5.1.5 to $Y=X^0, Z=\pi_0(X)$, so $\epsilon=\pi_0(i)$,
$\widetilde{\epsilon} = i^0$, and we keep the notations $j, \widetilde{j}$ for the
other two morphisms. We take $\xi = [\uH^\b(E^\b)]$. Then $j_*\xi=[E^\b]\in K_G(X^0)$.
Because $E^\b$ is, in particular, a complex of vector bundles on $X^0$, taking Tor's of
$\uH^\b(E^\b)$ with $\cO_{X^{0G}}$, as in (5.1.4), gives the same element of K-theory
as just tensoring $E^\b$ with $\cO_{X^{0G}}$, i.e., forming the restriction
$E^\b\big|_{X^{0G}}$. In other words,
$$\pi_0(i)^!\,\, [\uH^\b(E^\b)] = \left[ \uH^\b\left( E^\b\big|_{X^{0G}}\right)\right].
\leqno (5.3.2)$$
Note further that $\cal N$, being the normal bundle of $X^{0G}$ in $X^0$, is the same as
$\nu^0$, so ${\cal N}\big|_{\pi_0(X^G)}=\overline{\nu}^0$. 
So Theorem 5.1.5 gives
$$[\uH^\b(E^\b)] = \pi_0(i)_*\left( { \left[ \uH^\b\left( E^\b\big|_{X^{0G}}\right)\right]\over
\left[ \Lambda^\b \left( \overline{\nu}^0{}^*\right)\right]}\right).
\leqno (5.3.3)$$
To prove Theorem 5.3.1 it is enough therefore to prove the following equality
in $K_G(\pi_0(X^G))$:
$$[\Lambda^\b(\overline{\nu}^1{}'{}^*)]\cdot  \left[ \uH^\b\left(E^\b\big|_{X^{0G}}\right)\right] = 
[\uH^\b(i^* E^\b)]. \leqno (5.3.4)$$
Let $I^\b\subset \cO_X^\b$ be the dg-ideal of $X^G$, so $I^0\subset \cO_{X^0}$ is the ideal
of $X^{0G}$. Then we have
$$E^\b\big|_{X^{0G}} = E^\b/I^0 E^\b, \quad i^*E^\b = E^\b/I^\b E^\b = 
\left(E^\b\big|_{X^{0G}}\right) \bigl/ I^{\leq -1}\cdot \left( E^\b\big|_{X^{0G}}\right).
\leqno (5.3.5)$$
Further, the usual  interpretation of the conormal bundle holds in the dg-situation as well:
$I^\b /(I^\b)^2 = \nu^*$. Therefore
$I^{\leq -1}/(I^{\leq -1})^2 \cdot I^0 = (\nu^*)^{\leq -1}$, and we deduce for the $I^{\leq -1}$-adic
filtration:
$$(I^{\leq -1})^d\cdot  \left( E^\b\big|_{X^{0G}}\right)\biggl/ (I^{\leq -1})^{d+1}\cdot 
 \left( E^\b\big|_{X^{0G}}\right) =  i^*(E^\b) \otimes_{\cO_{X^G}^\b} S^d ((\nu^*)^{\leq -1}).
\leqno (5.3.6)$$
Notice that Corollary 4.4.10 is applicable to equivariant K-groups as well. Applying it to the
dg-variety $X^G$, we get
$$\bigl[ \uH^\b\bigl(i^*E^\b \otimes_{\cO^\b_{X^G}} S^d((\nu^*)^{\leq -1})\bigr)\bigr] = 
[S^d((\overline{\nu}^*)^{\leq -1})] \cdot [\uH^\b(i^*E^\b)]
= [S^d(\overline{\nu}^1{}'{}^*)]\cdot [\uH^\b(i^*E)],\leqno (5.3.7)$$
where the last equality follows from the quasiisomorphism of  $(\overline{\nu}^*)^{\leq -1}$
with $\overline{\nu}^1{}'{}^*$. 

Now, at the formal level, if we replace $E^\b\big|_{X^{0G}}$ by the (infinite) sum of the quotients
of the $I^{\leq -1}$-adic filtration, given by (5.3.6), we get the sum of the classes of the
symmetric
powers of $\overline{\nu}^1{}'{}^*$ which is (formally) inverse to the class of the
exterior algebra in (5.3.4). This can be made rigorous by performing the deformation to
the normal cone to $X^G$ in $X$, i.e., by considering the $I^\b$-adic filtration
in $\cO^\b_X$ and its associated graded sheaf of algebras $\gr_I\cO^\b_X$. 
Its spectrum is ${\cal N}_{X^G/X}$, the (dg)-normal bundle to $X^G$ in $X$ considered as a dg-manifold. 
Let us denote it $\widehat{X}$. 
Note that its underlying scheme $\widehat{X}^0$ is ${\cal N}_{X^{0G}/X^0}$, the normal
 bundle to $X^{0G}$
in $X^0$. At the same time $\widehat{X}^G = X^G$. Let $\widehat{i}: X^G\to \widehat{X}$ 
be the embedding. 
 Taking the $I$-adic filtration in $E^\b$, we have that $\gr_IE^\b$ is a module over
$\gr_I\cO^\b_X$ and thus gives a dg-vector bundle $\widehat{E}^\bullet$
on $\widehat{X}$. 
As in (3.3.3-4), the argument with a spectral sequence of coherent sheaves
on ${\cal N}_{X^{0G}/X^0}$, converging for Noetherian reasons, gives that
$$\bigl[\uH^\b(E^\b\big|_{X^{0G}}\bigr)\bigr] = \bigl[\uH^\b(\widehat{E}^\b\big|_{X^{0G}}\bigr)\bigr],
\quad [\uH^\b(i^*E)] = [\uH^\b(\widehat{i}^*\widehat E)].\leqno (5.3.8)$$ 
So we can and will  assume in proving (5.3.4), that $X=\widehat{X}$ coincides with
the normal bundle to the fixed point locus. 
In this case, the $I^\bullet$-adic filtration comes from a grading, so 
$$E^\b\big|_{X^{0G}} \quad = \quad \bigoplus_{d=0}^\infty \,\,\, (i^*E)\otimes_{\cO^\b_{X^G}}S^d((\nu^*)^{\leq -1}),$$
and the LHS of (5.3.4) becomes, by Corollary 4.4.10,
$$\biggl[ \uH^\b\biggl(\Lambda^\b((\nu^*)^{\leq -1})\otimes_{\cO^\b_{X^G}} S^\b ((\nu^*)^{\leq -1})\otimes
_{\cO^\b_{X^G}} i^*E\biggr)\biggr].\leqno (5.3.9)$$
Let $d$ be the differential in the triple tensor product of
complexes in (5.3.9). We can add to $d$ another summand $\delta$, the Koszul differential on $\Lambda^\b\otimes S^\b$
tensored with the identity  on the third factor, and we can arrange the tensor product into  a double complex.
The cohomology with respect to $\delta$ is then $i^*E$, so $\uH^\b_d(\uH^\b_{\delta}) = \uH^\b(i^*E)$,
and a spectral sequence argument shows that its class in $K_G(\pi_0(X^G))$ is the same
as the class of $\uH^\b_{d+\delta}$. On the other hand, the class of $\uH^\b_{d+\delta}$ is
equal to that of $\uH^\b_d$, as we see from the other spectral sequence corresponding to
the double complex. This proves the equality (5.3.4) and Theorem 5.3.1.

\vfill\eject

\baselineskip =12pt

\centerline {\bf References.}

\vskip 1cm

\noindent [ABW] K. Akin, D.A.  Buchsbaum, J. Weyman, 
{\it Schur functors and Schur complexes}, 
Adv. in Math. {\bf 44} (1982),  207-278.

\vskip .1cm

\noindent [B] K. Behrend, {\it Differential Graded Schemes I: Perfect Resolving Algebras},
preprint \hfill\break  math.AG/0212225.

\vskip .1cm

\noindent [BF] K. Behrend, B. Fantechi, {\it The intrinsic normal cone},
 Invent. Math. {\bf 128} (1997),  45-88.

\vskip .1cm

\noindent [CK1] I.  Ciocan-Fontanine, M.  Kapranov, {\it Derived Quot schemes,}
 Ann. Sci. \'Ecole Norm. Sup. (4) {\bf 34} (2001),  403-440. 

\vskip .1cm

\noindent [CK2]  I.  Ciocan-Fontanine, M.  Kapranov, {\it Derived Hilbert schemes,}
J. Amer. Math. Soc. {\bf 15} (2002),  787-815. 

\vskip .1cm

\noindent [F] W. Fulton, {\it Intersection Theory}, Springer-Verlag, 1984.

\vskip .1cm

\noindent [GP] T.  Graber, R.  Pandharipande,
{\it  Localization of virtual classes},  Invent. Math. {\bf 135} (1999),  487-518.

\vskip .1cm

\noindent [K] M.  Kontsevich, {\it  Enumeration of rational curves via torus actions}, 
in:  ``The Moduli Space of Curves'' (Texel Island, 1994), 335-368, Progr. Math., 
{\bf 129}, Birkh\"auser Boston, Boston, 1995. 

\vskip .1cm

\noindent [LT] J. Li and G. Tian, {\it Virtual moduli cycles
and Gromov-Witten invariants of algebraic varieties}, 
J. Amer. Math. Soc. {\bf 11} (1998),  119-174.

\vskip .1cm

\noindent [M] I. G. Macdonald, {\it Symmetric Functions and Hall polynomials},
Cambridge Univ. Press, 1986. 

\vskip .1cm

\noindent [T] R.W. Thomason, {\it Une formule de Lefschetz en K-th\'eorie
\'equivariante alg\'ebrique}, Duke Math. J. {\bf 68} (1992), 447-462. 

\vskip .1cm

\noindent [R] D.C. Ravenel, {\it Complex Cobordism and Stable Homotopy Groups of Spheres},
Academic Press, 1986. 

\vskip .1cm

\noindent [S] R. E. Stong, {\it Relations among characteristic numbers I,}
Topology, {\bf 4} (1965), 267-281.

\vskip 2cm

\noindent I.C.-F.: Department of Mathematics, University of Minnesota, 
127 Vincent Hall, 206 Church St. S.E., Minneapolis, MN 55455 USA, email:
$<$ciocan@math.umn.edu$>$

\vskip .3cm

\noindent M.K.: Department of Mathematics, Yale University, 10 Hillhouse Avenue,
New Haven, CT 06520 USA, email: $<$mikhail.kapranov@yale.edu$>$

\end